\documentclass[%
 reprint,
 amsmath,amssymb,
 aps,
]{revtex4-1}

\usepackage{graphicx}
\usepackage{dcolumn}
\usepackage{yhmath}
\usepackage{mathdots}
\usepackage{MnSymbol}
\usepackage{bm}
\usepackage{listings}
\usepackage{color}
\usepackage{subcaption}
\usepackage[toc,page]{appendix}

\usepackage{amsmath}
\usepackage{algorithm}
\usepackage[noend]{algpseudocode}

\makeatletter
\def\BState{\State\hskip-\ALG@thistlm}
\makeatother

\definecolor{dkgreen}{rgb}{0,0.6,0}
\definecolor{gray}{rgb}{0.5,0.5,0.5}
\definecolor{mauve}{rgb}{0.58,0,0.82}

\begin{document}

\preprint{APS/123-QED}

\title{Cleaning the correlation matrix with a denoising autoencoder}

\author{Soufiane Hayou}
 \altaffiliation[Also at ]{CMAP, Ecole Polytechnique (Paris).}
\affiliation{%
Quantitative Research, Bloomberg LP.\\
soufiane.hayou@polytechnique.edu
}%

\date{\today}

\begin{abstract}
In this paper, we use an adjusted autoencoder to estimate the true eigenvalues of the population correlation matrix from a noisy sample correlation matrix when the number of samples is small. We show that the model outperforms the Rotational Invariant Estimator (Bouchaud and Bun \cite{rmt}) which is the optimal estimator in the sample eigenvectors basis. 

\end{abstract}

\pacs{Valid PACS appear here}
\maketitle


\section{Introduction}

Correlation matrices have been used for decades to describe the joint dynamics of different variables. Although the correlation is just a measure of the 'linearity' between two variables, it is still very informative in many real-life problems : signal processing, statistical physics, portfolio optimization ... \\

It is well known that when we have a large number of data (compared to the dimension of the variables), we can accurately estimate the correlation matrix using the sample covariance matrix :\\
Let $N$ be the dimension, $T$ the number of data available and $(X_i)_{1\leq i \leq T}$ the observations. The sample covariance matrix is defined by :
\begin{equation*}
S = \frac{1}{T} \sum_{i=1}^T X_i X^t_i
\end{equation*}
where $X^t$ is the transpose of $X$. The empirical correlation is then calculated by dividing each component by $\sqrt{var(X_i) var(X_j)}$.\\

However, the sample covariance matrix is not a good estimator when $T$ is small (not much larger than $N$). Actually, for a large class of matrices (see \cite{soufiane}), we empirically observe that the sample covariance matrix tends to overestimate (underestimate) the large (small) eigenvalues. From now on, we only study the correlation matrix since the estimation of the variance is relatively robust using only the usual estimators.\\ 

There is an extensive amount of research papers on this topic. We find (in general) three main categories of approaches : Shrinkage to a target (e.g. Ledoit and Wolf \cite{ledoit}), Random Matrix Theory (N. El Karoui \cite{elkaroui}, Bouchaud and Bun \cite{rmt}) and Optimization under constraints (e.g. under a constraint on the condition number like in \cite{condition}). In this paper, we use a different approach based on machine learning and inspired from the results of Bouchaud and Bun on the estimation of the eigenvalues using a RIE (Rotational Invariant Estimator) and Random Matrix Theory.\\

\textbf{Contribution of this paper }: we show that an adjusted autoencoder (autoencoder with a noise level input, see section II) can outperform the RIE estimator in the case where $T$ is not vary large compared to $N$. The most important part is to have an exhaustive training dataset (exhaustive in terms of correlation matrices distribution). We use different simulation methods for this purpose.\\

In all the paper, S is the sample correlation matrix, C the population matrix, $\lambda_1 \geq \lambda_2 \geq ... \geq \lambda_N$ the eigenvalues of S and $u_1, u_2,...  ,u_N$ the corresponding eigenvectors, $c_1 \geq c_2 \geq ... \geq c_N$ the eigenvalues of C and $v_1, v_2,...  ,v_N$ the corresponding eigenvectors, and $q = \frac{N}{T} < 1$ is the ratio of the dimension over the number of samples.\\ 

In \cite{rmt}, authors show that the optimal (oracle) estimator of the true eigenvalues (in the sample eigenvectors basis) is given by an explicit formula in the large dimensional limit (see Appendix I). This formula is expressed using only the sample eigenvalues and the ratio $q = \frac{N}{T}$. However, this formula is supposed to work only when $N \rightarrow \infty$ (with $q$ fixed), and is the best estimator only in the sample eigenvectors space. However, we start from the observation that the only parameters needed in this estimator are the sample eigenvalues and the 'noise level' parameter $q = \frac{N}{T}$. The idea is to train an autoencoder to learn the mapping between the sample eigenvalues and the true eigenvalues where we feed also the parameter $q$ to the model.\\

In section II, we show the basic denoising autoencoder and the model we call the 'adjusted' autoencoder where we add a noise parameter to the input. We show also the numerical results (we use $L_2$ norm to estimate the accuracy). In section III, we show how to generate training data.

\newpage

\section{Denoising Autoencoder}
Let us begin by recalling the definition of the basic autoencoder (the one used in \cite{bengio}). In the following, X and Y are two random variables with joint probability distribution p(X,Y).\\

An autoencoder is a neural network where the input and the output have the same dimension. It is usually used to learn a representation of the original input for the purpose of dimensionality reduction. In our case, we use it for a different purpose (see below).

\subsection{Traditional Autoencoder}
The traditional autoencoder takes a vector x of dimension d as input and maps it to a hidden representation y (of different dimension) which is also mapped to an ouput z of dimension d :
\begin{figure*}
\begin{subfigure}{.5\textwidth}
  \centering
  \includegraphics[width=1\linewidth]{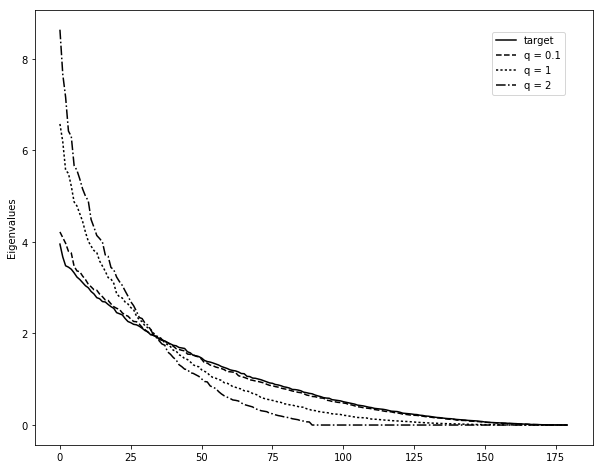}
  \caption{Example of exponentially decaying spectrum}
  \label{fig:sfig1}
\end{subfigure}%
\begin{subfigure}{.5\textwidth}
  \centering
  \includegraphics[width=1\linewidth]{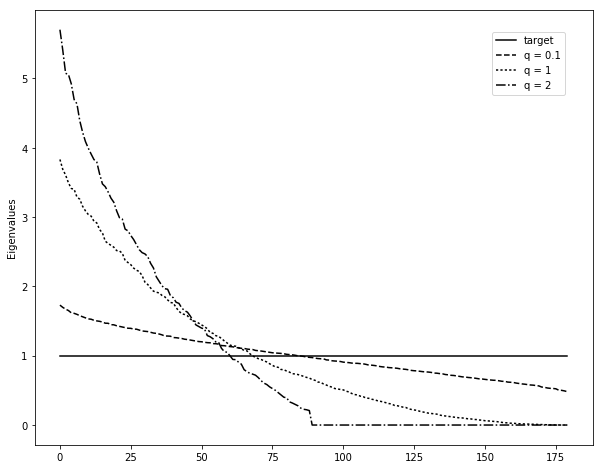}
  \caption{Flat spectrum}
  \label{fig:sfig2}
\end{subfigure}
\caption{Sample spectrum for different values of $q=\frac{N}{T}$}
\label{fig:fig}
\end{figure*}

\begin{center}
$y = a(Wx + b)$\\
$z = a(W'y+b')$
\end{center}

where $\theta = [W,b]$ and $\theta' = [W',b']$ are the model weights.\\

In this setting, $z$ can be viewed as a parameter of a distribution $p(X|Z=z)$ that generates X with hight probability. The error is defined then by : 
\begin{equation*}
L(x,z) \propto - log(p(x|z))
\end{equation*}\\

For real-valued variables, one example is : $X|z \sim N(z,\sigma^2 I)$ . This yields to :
\begin{equation}
L(x,z) = C ||x-z||^2
\end{equation}

To train the autoencoder, we minimise the Loss function with respect to the weights : 
\begin{equation}
\arg \min_{\theta, \theta'} \mathbb{E}_{X} [L(X,Z(X))]
\end{equation}

One natural hypothesis (constraint) on the traditional autoencoder is to assume that $W'=W^t$ where $W^t$ is the transpose of $W$.\\

When the dimension of the hidden layer is smaller than the input dimension, the autoencoder can be used as a data compressor (dimensionality reduction), and the hidden variable $z$ is the new representation of the input $x$. This is by far the most used application of autoencoders. A natural question is : why not using PCA instead ? the answer is that autoencoders are much more flexible, we can use different activation functions which adds non-linearities to the compression, whereas PCA can only use linear combinations of the variables.\\

\subsection{Denoising Autoencoder}
We have seen that the (trained) traditional autoencoder is a mapping of the identity with hidden representations of the input. Now what happens in the case where the input is corrupted (noisy observations) ? we would like to have a better representation of the input (cleaning) and that is exactly what a denoising autoencoder is trained to do. Of course the 'cleaning' accuracy will depend on the distribution of the noise around the 'true' input, and different noise distributions will lead to different 'cleaned' inputs.\\

To train the denoising auto-encoder, one usually uses a mapping distribution $x' \sim q(x'|x)$ and use it as the input to the model and the true input $x$ as the output :\\

- Map $x$ to $x' \sim q(x'|x)$\\

- Train the auto-encoder with inputs $x'$ and outputs $x$.\\

There is many types of distribution that can be used to corrupt the input. We cite the ones used in \cite{vincent} : \\

- Additive isotropic Gaussian noise : $x'|x \sim N(x, \sigma^2 I)$\\

- Masking noise : some elements of $x$ are forced to be 0.\\

- Salt-and-paper noise : some elements of $x$ are set to their maximum/minimum value possible.\\

The Additive gaussian noise is a common noise model for real-valued variables, whereas the Salt-and-paper noise model is a natural choice for binary (or almost binary) variables.
\begin{figure*}
  \includegraphics[width=0.8\linewidth]{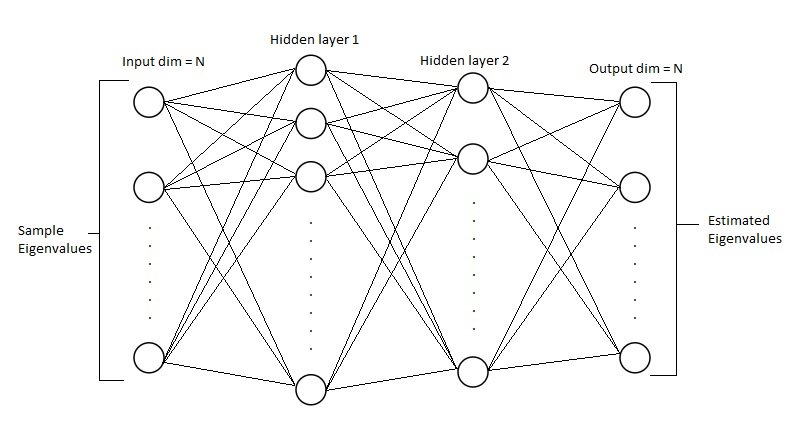}
  \caption{Denoising autoencoder}
\end{figure*}

In \cite{vincent}, authors use a denoising autoencoders to clean the input of the layers (adding a denoising autoencoder before the layer). The model is called Stacked Denoising Autoencoders,

\section{The model}
We recall that $q = \frac{N}{T}$ is the ratio of the dimension over the number of samples. We empirically observe that this parameter encodes the noise level around the true eigenvalues (the more data we have the better the estimation).\\

\subsection{Noise level}
The sample correlation matrix is a noisy version of the true correlation matrix. Figure 1 shows an example of the noise level for a correlation matrix with dimension $N = 180$. The blue line is the true spectrum, while the other lines are the estimates of the spectrum from different numbers of samples ($T = 1800$, $T=360$ and $T=180$).

This shows that a 'good' estimator of the true eigenvalues should use $q$ as a parameter(which we find in Bouchaud's formula of the optimal estimator, see Appendix 1).\\

Now let us go back to the original problem. We want to clean the spectrum of the sample correlation matrix in order to approximate the true spectrum. That means this problem reduces to a Denoising problem, and the idea of using a Denoising Autoencoder becomes natural. However, the difference here is that we don't know the distribution of the noise around the true spectrum and thus we cannot directly generate $x'$ (the corrupted version of $x$, note that x here is the vector of eigenvalues, which means the autoencoders has an input/output dimension of $N$). We solve this problem indirectly by simulating data from the correlation matrix and compute the corrupted spectrum from the sample correlation matrix calculated with the data. 

\subsection{Cleaning with a denoising autoencoder}
The purpose of this section is to show that the usual denoising autoencoder will perform better than RIE for a unique value of q (the dimension of number of samples are fixed), and will perform poorly for values of $q$. The generation of training data will be discussed in the next section.\\
Figure 2 shows the autoencoder we use for this purpose. We use a autoencoder with two hidden layers (we use two hidden layers instead of one because depth gives more approximation power, see \cite{ronen}). Note that when the autoencoder is used for a compression purpose, the hidden layers have usually less neurons than the input dimension, which is not the case here, since we are trying to learn the true curvature of the spectrum (and not a compression of the data). In our example, we have $N = 180$, $n_1 = 300$ (number of neurons in the first hidden layer) and $n_2= 200$ (number of neurons in the second layer), and to avoid over-fitting, we add a 'Dropout' on the second layer. \\
To show the impact of q, we train the model with a dataset generated with a fixed $q$ (same number of samples used to calculate the sample correlation matrix).\\

Here the number of samples generated from the true correlation matrix is fixed to $T=200$ ($q = 180/200 = 0.9$). In the example below, we compare the autoencoder output with Bouchaud's estimation and the sample estimation (the comparison is done with the $L_2$ norm).

\begin{figure}
  \includegraphics[width=1.\linewidth]{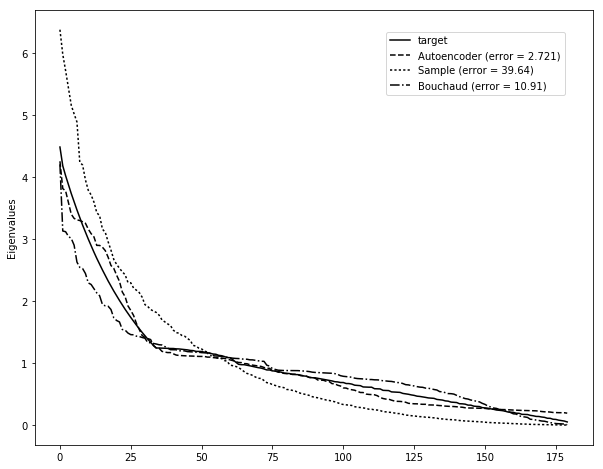}
  \caption{Example of estimating the true eigenvalues using an autoencoder with $T = 190$ (near $q$ used for the training)}
\end{figure}

In Figure 3, the autoencoder outperforms the other methods in terms of the $L_2$ norm. Now we test the model with T = 400.
\begin{figure}[!htb]
  \includegraphics[width=1.\linewidth]{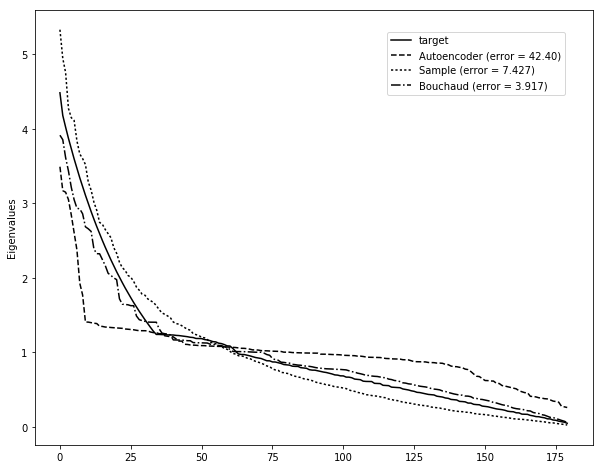}
  \caption{Example of estimating the true eigenvalues using an autoencoder for $T = 400$ (different $q$ from the one used for the training)}
\end{figure}

Figure 4 shows that the autoencoder performs poorly when we use a value of $q$ (we change $T$) different from the one used in the training step. Bouchaud's estimator is the best estimator in this case.\\

\begin{figure}[!htb]
  \includegraphics[width=1.\linewidth]{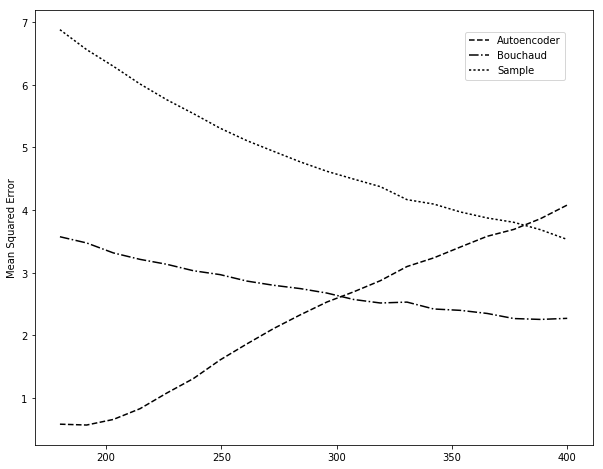}
  \caption{MSE for $T$ between 180 and 400 ($q$ between 1 and 0.45)}
\end{figure}

In Figure 5, we calculate the mean squared error for a range of different values of $q$ (different values of $T$). It shows that the autoencoder performs poorly when $q$ is far from the value of $q$ used for the training.  

\begin{figure*}
  \includegraphics[width=0.8\linewidth]{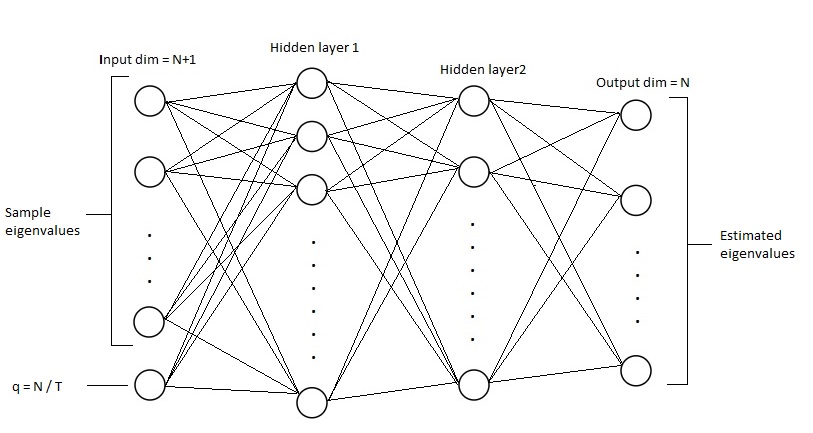}
  \caption{The adjusted autoencoder}
\end{figure*}

\begin{figure*}
  \includegraphics[width=0.5\linewidth]{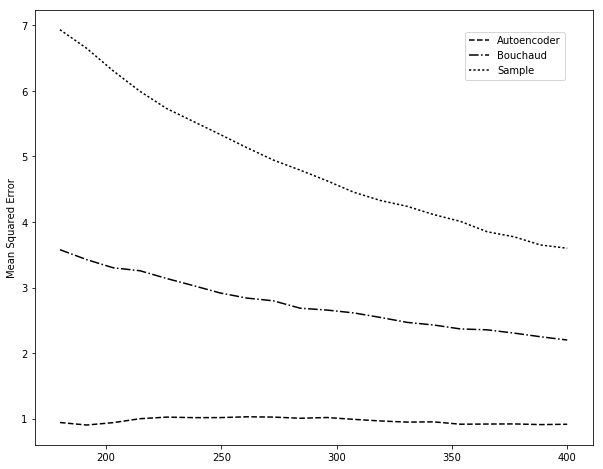}
  \caption{MSE for $T$ between 180 and 400 ($q$ between 1 and 0.45)}
\end{figure*}

\subsection{The adjusted Autoencoder}

The idea is to add the noise parameter $q$ as an input to the autoencoder. That means instead of having an input of dimension $N$ we will have an input of dimension $N+1$ and an output of dimension $N$. We call the new model the 'adjusted' autoencoder. We use a fully connected autoencoder to which we add another input ($q = N/T$). \\

The sample eigenvalues and the true eigenvalues are sorted increasingly before they are fed to the model (the intuition behind this is to distinguish the spectrum from the parameter $q$). This makes the model learn the function that maps the curvature of the sample spectrum to the curvature of the true spectrum. We add also a Dropout on the second layer (with dropout probability 25$\%$). We choose $N = 180$, $n_1 = 300$ (number of neurons in the first hidden layer) and $n_2= 200$ (number of neurons in the second layer), Figure 6 shows the final model.\\

To see the impact of adding the parameter $q$ to the training, we show in Figure 7 the mean squared error for different values of $q$ (different values of $T$). We use out-of-sample data to do the calculations. It shows that the adjusted autoencoder performs better than the RIE for a wide range of values of $T$ (the model was trained for different values of T, see next section).\\

In the next page, we show different examples of the estimation of the eigenvalues using the adjusted autoencoder.
\begin{figure*}
\begin{subfigure}{.5\textwidth}
  \centering
  \includegraphics[width=1\linewidth]{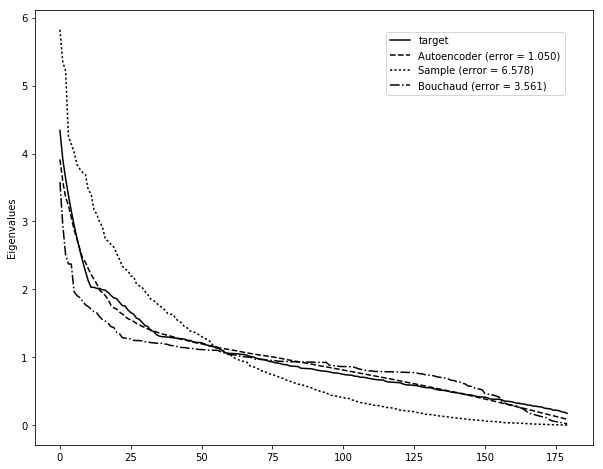}
  \caption{Exponentially decaying spectrum}
  \label{fig:sfig1}
\end{subfigure}%
\begin{subfigure}{.5\textwidth}
  \centering
  \includegraphics[width=1\linewidth]{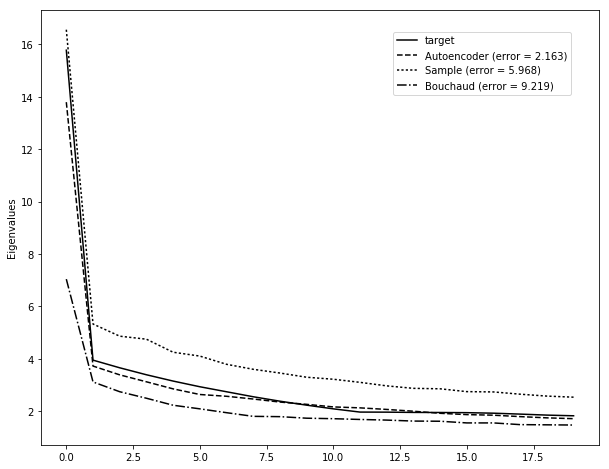}
  \caption{Spiked spectrum}
  \label{fig:sfig2}
\end{subfigure}
\begin{subfigure}{.5\textwidth}
  \centering
  \includegraphics[width=1\linewidth]{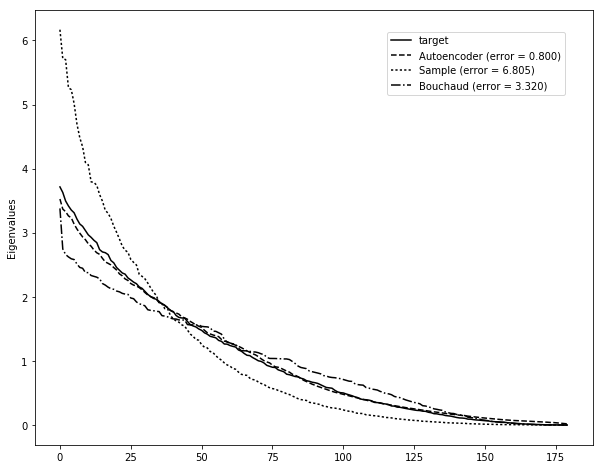}
  \caption{Slow decaying spectrum}
  \label{fig:sfig2}
\end{subfigure}%
\begin{subfigure}{.5\textwidth}
  \centering
  \includegraphics[width=1\linewidth]{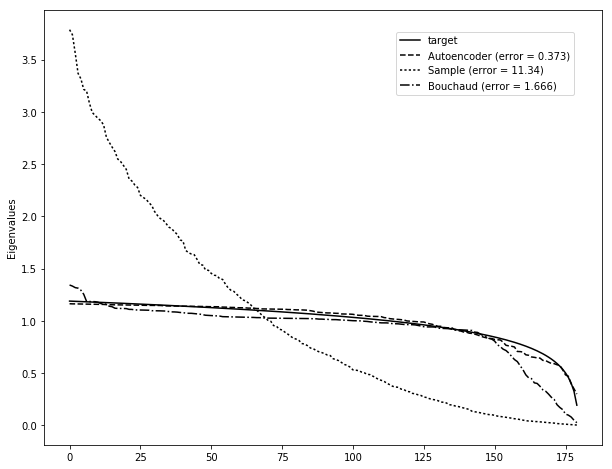}
  \caption{Concave spectrum}
  \label{fig:sfig2}
\end{subfigure}
\begin{subfigure}{.5\textwidth}
  \centering
  \includegraphics[width=1\linewidth]{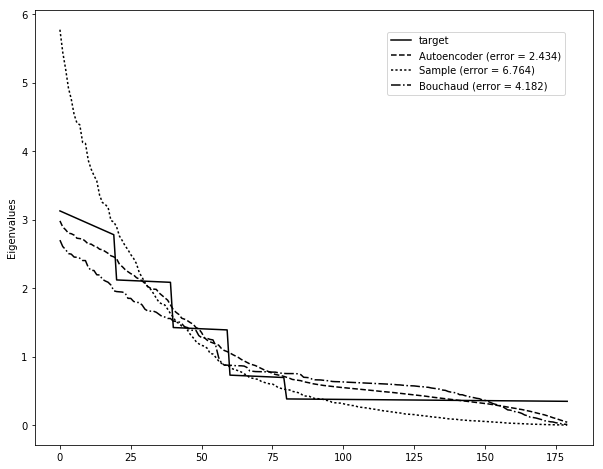}
  \caption{Extreme scenario 1}
  \label{fig:sfig2}
\end{subfigure}%
\begin{subfigure}{.5\textwidth}
  \centering
  \includegraphics[width=1\linewidth]{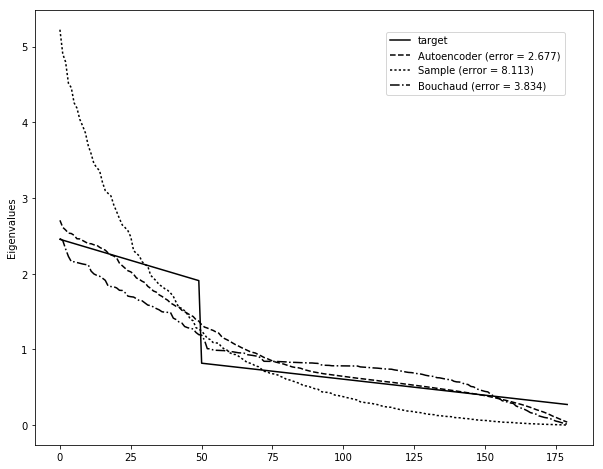}
  \caption{Extreme scenario 1}
  \label{fig:sfig2}
\end{subfigure}
\caption{Output of the adjusted autoencoder for different examples of spectrum}
\label{fig:fig}
\end{figure*}

\newpage

\section{Generating training data}
In this section, we present some methods for the generation of training data. \\
From a given correlation matrix (and since the correlation matrix does not depend on the variance of the variables), we generate variables with unit variances and take the sample covariance matrix (which is also the sample correlation matrix in this case). Note that we only need the spectrum (and not the correlation matrix itself) to generate the samples (see Appendix II). But since we are working on correlation matrices, we present some methods of generation of a random correlation matrix.  

\subsection{Generating a random correlation matrix with specified eigenvalues}
This generation method was first introduced by Davies and Higham \cite{higham}. In the following we give a brief description of the algorithm.\\

The algorithm is divided into 3 steps : generating a random orthogonal matrix $\rightarrow$ generating a random matrix with specified eigenvalues $\rightarrow$ applying Given's rotations to transform the previous matrix into a correlation matrix. \\

The generation of a random orthogonal matrix was the subject of an extensive amount of research papers, many of them propose a simulation method. In this paper we use Stewart's algorithm \cite{stewart}  to generate a random orthogonal matrix.

\subsubsection*{Generating a random matrix with specified eigenvalues}
The next step is to use the orthogonal matrix (previously generated) to construct a matrix that has some specified eigenvalues. This is straightforward using the following formula : \\

Let $a_1, a_2, ..., a_n$ the eigenvalues and Q an orthogonal matrix. Then :\\
\begin{equation}
M = Q^{t} Diag(a_1, a_2, ..., a_n) Q
\end{equation}
where $Q^t$ is the transpose of $Q$, is a matrix with eigenvalues $a_1, a_2, ..., a_n$.

\subsubsection*{Given's rotations}
Now that we have a matrix M with eigenvalues $a_1, a_2, ..., a_n$ (note that we should have $\sum_{i=1}^n a_i = n$ so that we can construct a correlation matrix with these eigenvalues), we can use Given's rotations to have 1's on the diagonal. Given's rotations are orthogonal transformations, ans thus the resulting matrix will have the same eigenvalues. In Appendix II, we present a short algorithm to do the Given's rotation on M in the $(i,j)$ position. \\

\subsubsection*{final algorithm}
Now that we have the Given's rotation tool, we can apply it on the matrix M (on the diagonal until we have only 1's), the resulting matrix will be a correlation matrix. In Annex '''', the function $generate\_corr\_wse$ ('generate correlation matrix with specified eigenvalues') returns a correlation matrix.\\

Note that the randomness in this algorithm comes from the randomness of the orthogonal matrix. So, different simulation methods of the orthogonal matrix may result in different density distributions of the resulting correlation matrix.\\

Let's prove that the output is a correlation matrix. For this purpose, we prove a more general lemma :\\

{\it If A is a positive symmetric matrix with 1's on the diagonal, then all the coefficients are between -1 and 1}\\

{\bf Proof :}\\
Let $(e_i)_{1 \leq i \leq n}$ be the usual basis of $\mathbb{R}^n$, and $A= (a_{ij})$ a positive symmetric matrix with 1's on the diagonal.\\
We have $a_{ij} = e^T_i A e_j$ for all $i,j$, and by the positiveness of A we have $(e_i - e_j)^T A (e_i - e_j) \geq 0$. We expand the left-hand formula :\\

$(e_i - e_j)^T A (e_i - e_j) = e_i^T A e_i + e_j^T A e_j - 2 e^T_i A e_j$\\

and using the fact that $e_i^T A e_i = e_j^T A e_j = 1$, we deduce that : \\
\begin{center}
$e^T_i A e_j \leq 1$.\\
\end{center}
We use the same idea with $(e_i + e_j)^T A (e_i + e_j)$ and we find that $e^T_i A e_j \geq -1$.\\

Conclusion : $|a_{ij}| = |e^T_i A e_j| \leq 1$ for all $i, j$.

\subsection{Other simulation methods}
In order to have an exhaustive training dataset, we added other simulation methods. Most of these methods were exposed in G. Marsaglia and I. Olkin \cite{olkin} and J. Hardin, S.R. Garcia, D. Golan \cite{hardin}. Here we present some of them.

\subsubsection*{Random correlation matrix of the form $AA^t$}

{\bf Lemma :} $AA^t$ is a random correlation matrix if and only if the rows of $A$ are random vectors on the unit sphere.\\

So, an easy way to generate a random correlation matrix is to generate a random matrix A (component-wise) and normalize the rows to be in the unit sphere. $AA^t$ is a random correlation matrix. Note that the distribution of the resulting correlation matrix depends on the initial distribution of the components of A. 

\subsubsection*{Generating random correlation matrices from constant correlation blocks}
In order to generate a random correlation matrix with dimension n, we use the following recipe.\\

Let K be the number of blocks (groups of constant correlation). For $k = 1,2, ..., K$, let $g_k$ be the size of the kth block (we should have $\sum_{k=1}^K g_k = n$), $\rho_k$ such that $0\leq \rho_k < 1$ the correlation of the group. We define $\rho_{min} = min{\rho_1, \rho_2, ..., \rho_K}$ and $\rho_{max} = max{\rho_1, \rho_2, ..., \rho_K}$. Let $\delta$ be a number such that $0\leq \delta < \rho_{min}$ (correlation between group).\\

We define the matrices $\Sigma_k$ (constant correlation blocks) by :
\[
M=
  \begin{bmatrix}
    1 & \rho_k & ... & \rho_k \\
    \rho_k & 1 & ... & \rho_k \\
    . & . & \ddots & . \\
    \rho_k & \rho_k & ... & 1 \\
  \end{bmatrix}
\]

Let $\sigma$ be the matrix with blocks $\Sigma_K$ on the diagonal and zeros elsewhere (dimension of $\sigma$ is n), $\epsilon$ be a real number such that $0 \leq \epsilon < 1- \rho_{max}$, and $x_1, x_2, ..., x_n$ n randomly generated unit vectors with dimension n. The matrix $Corr$ defined by :

\[
  Corr_{i,j} = \left\{\def\arraystretch{1.2}%
  \begin{array}{@{}c@{\quad}l@{}}
    1 & \text{if i=j}\\
    \rho_k + \epsilon x_i^T x_j & \text{if i,j are in the kth group and i$\neq$j}\\
    \delta + \epsilon x_i^T x_j & \text{if i,j are in different groups}\\
  \end{array}\right.
\]
is a correlation matrix with the following upper bound on its condition number : 
\begin{equation}
\kappa(Corr) \leq \frac{n(1+\epsilon) + 1}{1 - \rho_{max} - \epsilon}
\end{equation}

\subsubsection*{Generating random correlation matrices from Toeplitz blocks}
We use the same notations of the previous method (K,k, $g_k$, n, $\rho_k$, $\Sigma$). The blocks have now the Toeplitz structure : 
\[
M=
  \begin{bmatrix}
    1 & \rho_k & \rho_k^2 &... & \rho_k^{g_k -1} \\
    \rho_k & 1 & \rho_k &... & \rho_k^{g_k -2} \\
    \rho_k & \rho_k^2 & 1 &... & \rho_k^{g_k -3} \\
    . & . & . &  \ddots & . \\
    \rho_k^{g_k-1} & \rho_k^{g_k-2} & ... & \rho_k & 1 \\
  \end{bmatrix}
\]

Let $\epsilon$ be a real number such that $0 < \epsilon < \frac{1 - \rho_{max}}{ 1 + \rho_{max}}$ and $x_1, x_2, ..., x_n$ n randomly generated unit vectors with dimension n. The matrix $Corr$ defined by : 
\begin{equation}
Corr = \Sigma + \epsilon (X^t X -I)
\end{equation}
where $X$ is the matrix with columns $x_i$ and $I$ is the identity matrix, is a correlation matrix, with the following upper-bound on the condition number :
\begin{equation}
\kappa(Corr) \leq \frac{\frac{1 + \rho_{max}}{1 - \rho_{max}} + (n - 1) \epsilon}{\frac{1 - \rho_{max}}{1 + \rho_{max}} - \epsilon}
\end{equation}\\

\subsection{Generating training data}
Now that we have the tools to generate random correlation matrices, we can use a combination of them to generate training data for the model. Note that we can directly generate training data using only the eigenvalues without using the correlation matrix (see Appendix II), and we don't use Given's rotation in this case. But since we generate the data just once, we use the previous algorithms.\\

It is clear that the first algorithm is exhaustive in the way that we can generate a correlation matrix with any given eigenvalues (that sum to $N$). To generate spectrums, we start by splitting the spectrum into 2 groups : Principal eigenvalues (biggest eigenvalues), and Other eigenvalues. We first generate a random number $p$ in $[0,1]$, p will be the percentage of variance explained by the Principal eigenvalues. Here is a sketch of the algorithm :\\

1. Generate uniformly in $[0,1]$ a real number $p$\\
2. Generate uniformly in $\{1,2,...,N\}$ an integer $l$\\
3. Generate uniformly in $[0,1]$ l numbers (Principal eigenvalues)\\
4. Scale the Principal eigenvalues so that the sum equals to p\\
5. Generate uniformly in $[0,1]$ N-l numbers (Other eigenvalues)\\
6. Scale the Other eigenvalues so that the sum equals to 1-p\\
7. Concatenate the two groups and scale the output so that the sum equals to N\\

Remark : since we don't know exactly the distribution of eigenvalues of a random correlation matrix, we add a combination of the other methods to the previous algorithm.

\cleardoublepage
\section*{Appendix I : Optimal Rotational Invariant Estimator}
In \cite{rmt}, the optimal RIE estimator (when $N \rightarrow \infty$) is given by :\\
\begin{equation*}
\xi^{ora}( \lambda) = \frac{1}{|1 - q + q \lambda \lim_{\nu \rightarrow 0^+ g_S(\lambda - i \nu)}|^2}
\end{equation*}

where,
\begin{equation*}
g_S(z) = \frac{1}{N}\textnormal{Tr}(z I_N - S)
\end{equation*}

is the Stieltjes transform of S.

\section*{Appendix II : Distribution of the sample eigenvalues}
{\bf Definition} : A $p \times p$ matrix M is said to have a Wishart distribution with covariance matrix $\Sigma$ and degrees of freedom n if $M = X^t X$ where $X \sim N_{n \times p}(\mu, \Sigma)$. We denote this by $M \sim W_p(n, \Sigma)$.\\

When $n \geq p$, the Wishart distribution has a density function given by : 
\begin{equation}
f(M) = \frac{2^{-np/2}}{\Gamma_p(n/2) (det(\Sigma))^{n/2}} etr(-\frac{1}{2} \Sigma^{-1} M) (det M)^{(n-p-1)/2}
\end{equation}
where $etr$ is the exponential of the trace, $\Gamma_p$ is the generalized gamma function.\\

When $X \sim N_{n \times p}(\mu, \Sigma)$, the sample covariance matrix $S = \frac{1}{n} X X^t$ has the Wishart distribution $W_p(n-1, \frac{1}{n}\Sigma)$.\\

{\bf Joint distribution of the eigenvalues}\\
Let $M \sim W_p(n,\Sigma)$ ($n >p$), then the joint distribution of the eigenvalues $l_1 \geq l_2 \geq ... \geq l_p$ is : 
\begin{equation}
\frac{\pi^{p^2/2} \times 2^{-np/2} (det \Sigma)^{-n/2}}{\Gamma_p(n/2) \Gamma_p(p/2)} \prod_{i=1}^p l_i^{(n-p-1)/2} \prod_{j>i}^p (l_i - l_j) \int_{O_p} etr(-\frac{1}{2} \Sigma^{-1} HLH^t) (dH)
\end{equation}

where the integral is over the orthogonal group $O_p$ with respect to the Haar measure (see \cite{Elizabeth}).\\
In general, the integral is hard to calculate, but in the case where $\Sigma = \lambda I$, we have : 
\begin{align*}
\int_{O_p} etr(-\frac{1}{2} \Sigma^{-1} HLH^t) (dH) &= \int_{O_p} etr(-\frac{1}{2 \lambda} HLH^t) (dH)\\
&= etr(-\frac{1}{2 \lambda}L) \int_{O_p} (dH)\\
&= exp(-\frac{1}{2 \lambda} \sum_{i=1}^p l_i)
\end{align*}

The Haar measure is invariant by rotation, that means for any orthogonal matrix $Q$, one has :
\begin{equation*}
d(QH) = dH
\end{equation*}

Using this and the fact that it exist an orthogonal matrix $Q$ such that $\Sigma^{-1} = Q D^{-1} Q^t$ where $D = diag(\lambda_1, \lambda_2, ..., \lambda_p)$, one can see that the previous distribution depends only on the eigenvalues of $\Sigma$.\\

\cleardoublepage

\section*{Appendix III : pythoon code for the generation of a random correlation matrix}

\subsection*{Given's rotation}

\begin{lstlisting}
def givens(M, i, j):
    G = M
    Mii, Mij, Mjj = M[i,i], M[i,j], M[j,j]
    t = (Mij + np.sqrt(Mij**2 - (Mii-1)*(Mjj-1))) / (Mjj - 1)
    c = 1. / np.sqrt(1+t**2)
    s  = c*t
    Mi, Mj = M[i], M[j]
    G[i], G[j] = c*Mi - s*Mj, s*Mi + c*Mj
    Mi, Mj = G[:,i], G[:,j]
    G[:,i], G[:,j] = c*Mi - s*Mj, s*Mi + c*Mj
    return G
\end{lstlisting}
\newpage
\subsection*{Generating random correlation matrices}
\begin{lstlisting}
def generate_corr_wse(eigs):
    n = len(eigs)
    eigen = n * eigs / np.sum(eigs)
    corr = generate_wge(eigen)
    precision = 0.1
    converg = 0
    i = 0
    while(not converg):
        vec = np.diagonal(corr)
        if np.sum(abs(vec-1)>precision)==0:
            converg = 1
        else:
            bigger = np.arange(len(vec))[(vec>1)]
            smaller = np.arange(len(vec))[(vec<1)]
            i,j = smaller[0], bigger[-1]
            if i>j:
                i,j = bigger[0], smaller[-1]
            corr = givens(corr, i, j)
            corr[i,i]=1
    return corr
\end{lstlisting}

\end{document}